\documentclass[a4paper]{article}
\usepackage{fullpage}
\usepackage{amsmath, amsthm, amsfonts, amssymb, enumerate}

\newtheorem{thm}{Theorem}[section]
\newtheorem{cor}[thm]{Corollary}
\newtheorem{lem}[thm]{Lemma}
\newtheorem{prop}[thm]{Proposition}

\newtheorem{que}[thm]{Question}
\newtheorem{defn}[thm]{Definition}
\newenvironment{de}{\begin{defn} \rm}{\end{defn}}

\newcommand{\aut}{\operatorname{Aut}}
\newcommand{\supp}{\operatorname{supp}}
\newcommand{\stab}{\operatorname{Stab}}
\newcommand{\fix}{\operatorname{fix}}
\newcommand{\autq}{\aut(\Q, \leq\nobreak)}
\newcommand{\auto}{\aut(\Omega, \leq\nobreak)}
\newcommand{\sym}{\operatorname{Sym}}

\newcommand{\Sq}[1]{\left( #1 \right)_{n \in \N}}
\newcommand{\N}{\mathbb{N}}
\newcommand{\Z}{\mathbb{Z}}
\newcommand{\Q}{\mathbb{Q}}
\newcommand{\R}{\mathbb{R}}
\newcommand{\id}{\operatorname{id}}
\newcommand{\cf}{\operatorname{cf}}
\newcommand{\scf}{\operatorname{scf}}
\newcommand{\set}[2]{\{#1:#2\}}
\renewcommand{\to}{\longrightarrow}

\begin{document}

\author{J. Hyde, J. Jonu\v sas, J. D. Mitchell and Y. P\'eresse}
\title{Universal sequences for the order-automorphisms of the rationals}

\maketitle
\begin{abstract}
  In this paper, we consider the group $\autq$ of order-automorphisms of the
  rational numbers, proving a result analogous to a theorem of Galvin's for the
  symmetric group.  In an announcement, Kh\'elif states that every countable
  subset of $\autq$ is contained in an $N$-generated subgroup of $\autq$ for
  some fixed $N\in\N$.  We show that the least such $N$ is $2$. Moreover, for
  every countable subset of $\autq$, we show that every element can be given as a
  prescribed product of two generators without using their inverses.  More
  precisely, suppose that $a$ and $b$ freely generate the free semigroup
  $\{a,b\}^+$ consisting of the non-empty words over $a$ and $b$.  Then we show
  that there exists a sequence of words $w_1, w_2,\ldots$ over $\{a,b\}$ such
  that for every sequence $f_1, f_2, \ldots\in \autq$ there is a homomorphism
  $\phi:\{a,b\}^{+}\to \autq$ where $(w_i)\phi=f_i$ for every $i$. 
  
  As a corollary to the main theorem in this paper, we obtain a result of
  Droste and Holland showing that the strong cofinality of $\autq$ is
  uncountable, or equivalently that $\autq$ has uncountable cofinality and
  Bergman's property.

\end{abstract}

%%%%%%%%%%%%%%%%%%%%%%%%%%%%%
%%%%%%%%%%%%%%%%%%%%%%%%%%%%%

\section{Introduction}

In \cite{Galvin1995aa}, Galvin shows that every countable subset of the
symmetric group $\sym(\Omega)$, on any infinite set $\Omega$, is contained in a
$2$-generated subgroup. The orders of the two generators can be chosen to be
almost any values, and in particular, the orders of both of the generators can
be finite.   It follows that the elements of any countable subset of
$\sym(\Omega)$ can be obtained as compositions of just $2$ permutations without
the use of their inverses.  In other words, Galvin obtained the following
theorem. 

%%%%%%%%%%%%%%%%%%%%%%%%%%%%%

\begin{thm}[cf. Theorem 4.1 in Galvin \cite{Galvin1995aa}]
  \label{thm-galvin}
  Let $\Omega$ be an arbitrary infinite set. Then every countable subset of
  the symmetric group $\sym(\Omega)$ on $\Omega$ is contained in a $2$-generated
  subsemigroup of $\sym(\Omega)$.
\end{thm}

%%%%%%%%%%%%%%%%%%%%%%%%%%%%%

A bijection $f:\Q\to \Q$ is an \emph{order-automorphism} when: $x\leq y$ if
and only if $(x)f\leq (y)f$. We denote the group of 
order-automorphisms of $\Q$ by $\autq$. In this paper we prove an analogue of
Theorem \ref{thm-galvin} for $\autq$. In an announcement Kh\'elif
\cite[Theorem 7]{Khelif2006aa}, states that every countable subset of $\autq$ is
contained in an $N$-generated subgroup of $\autq$ for some fixed $N\in\N$.
However, there is no proof of Kh\'elif's assertion in \cite{Khelif2006aa}, nor
is the value of $N$ mentioned. We give a proof of Kh\'elif's assertion showing
that $N$ can, in fact, be $2$.

A \emph{transformation} of a set $\Omega$ is simply any function from $\Omega$
to itself.  Galvin was motivated by the following two theorems and a question of
Stan Wagon, who asked if ``transformation'' could be replaced by ``permutation''
in Theorem \ref{thm-sierpinski}. 

%%%%%%%%%%%%%%%%%%%%%%%%%%%%%

\begin{thm}[Theorem IV in Higman-Neumann-Neumann \cite{Higman1949aa}]
  Every countable group is embeddable in a 2-generated group.
\end{thm}

%%%%%%%%%%%%%%%%%%%%%%%%%%%%%

\begin{thm}[Sierpi\'nski \cite{Sierpinski1935aa} and Banach \cite{Banach1935aa}]
  \label{thm-sierpinski}
  Every countable set of transformations on an infinite set $\Omega$ is
  contained in a semigroup generated by two transformations of $\Omega$.
\end{thm}

%%%%%%%%%%%%%%%%%%%%%%%%%%%%%

Analogues of Sierpi\'nski's theorem have been found for several further classes
of groups and semigroups; see the introduction to \cite{Mitchell2011ab} for more
details.  Perhaps most relevant for our purposes is that Galvin's proof can be
adapted to show that if $G$ is the group of homeomorphisms of the Cantor space,
the rationals, or the irrationals, then any countable subset of $G$ is contained
in a $2$-generator subgroup.  It was shown by Calegari, Freedman, and de
Cornulier \cite{Calegari2006aa} that the homeomorphisms of the euclidean
$m$-sphere have the property that every countable subset is contained in a
$N$-generated subgroup for some $N\in \N$ (the specific value of $N$ is not
given in \cite{Calegari2006aa}). To our knowledge, these examples exhaust the
naturally arising non-finitely generated groups which are known to have the
property that every countable subset is contained in an $m$-generated subgroup,
for some fixed $m\in\N$.  This property is preserved under taking subgroups of
finite index, finite direct products, and arbitrary restricted wreath products,
which give rise to further examples.

%%%%%%%%%%%%%%%%%%%%%%%%%%%%%

The property investigated in this paper is somewhat stronger than the
property mentioned above and will be defined next. We write $A^+$ to denote the free semigroup freely generated by an alphabet $A$ and refer to a sequence of elements of $A^+$ as \emph{a sequence over $A$}.

\begin{de}
  Let $S$ be a semigroup, let $T$ be a subset of $S$,
  and let $A$ be an alphabet.  Then a sequence $w_1, w_2,\ldots$ over $A$ is \emph{universal
  for $T$ as a subset of $S$} if for any sequence $t_1, t_2, \ldots \in T$ there exists a
  homomorphism $\phi:A^+\to S$ such that $(w_i)\phi=t_i$ for all $i\in \N$. 
\end{de}

If the alphabet $A$ has $m$ elements, then we will refer to $\Sq{w_n}$ as an
\emph{$m$-letter universal sequence}.  If a sequence is universal for $S$ as a subset of
$S$, then we simply refer to this sequence as \emph{universal for $S$}.
Wherever it is possible to do so without ambiguity, we will also not refer
specifically to the alphabet $A$.  Of course, the analogous definition of
universal sequences for groups can be given using the free group.  However, we
will not use the analogous definition in this article.  

In an announcement, Kh\'elif \cite{Khelif2006aa} states that there is a finite
letter universal sequence for $\autq$. However, neither the number of letters in
the universal sequence, nor a proof of this statement, is given in
\cite{Khelif2006aa}.

In this paper, in the spirit of Galvin's Theorem, we will show that $\autq$ has
a $2$-letter universal sequence.  In other words, every element from an
arbitrary countable set can be given as a prescribed product of two generators
without using their inverses.  The main result of this paper is the following
theorem. \medskip

%%%%%%%%%%%%%%%%%%%%%%%%%%%%%

\noindent\textbf{Main Theorem} (cf. Theorem 7 in Kh\'elif \cite{Khelif2006aa}).
\textit{There is a $2$-letter universal sequence for $\autq$.}\medskip

%%%%%%%%%%%%%%%%%%%%%%%%%%%%%

We conclude the introduction by discussing some of the consequences of the Main
Theorem.  The \emph{cofinality} of a group $G$, denoted $\cf(G)$, is the
least cardinal $\lambda$ such that $G$ can be written as the union of a chain of
$\lambda$ proper subgroups. Macpherson and Neumann \cite{Macpherson1990aa}
showed that the symmetric group on a countably infinite set has uncountable
cofinality; Gourion \cite{Gourion1992aa} showed that $\cf(\autq)>\aleph_0$; Hodges
\emph{et al.} \cite{Hodges1993aa} and Thomas \cite{Thomas1996ab} proved the
analogous results for the automorphism group of the random graph and the
infinite dimensional linear groups, respectively. 

The \emph{strong cofinality} of a group $G$, denoted $\scf(G)$, is the least
cardinal $\lambda$ such that $G$ can be written as the union of a chain of
$\lambda$ proper subsets $H_i$ such that for all $i$ the following hold:
\begin{itemize}
  \item $H_i={H_i}^{-1}$;
  \item there exists $j\geq i$ with $H_iH_i\subseteq H_j$.
\end{itemize}
Droste and Holland \cite{Droste2005aa} showed that $\scf(\autq)>\aleph_0$.  A
group $G$ has \emph{Bergman's property} if for any generating set $X$ for $G$
with $X=X^{-1}$ and $1\in X$, there is $N\in \N$ such that $G=X^N$.  Bergman
\cite{Bergman2006aa} showed that the symmetric group $\sym(\Omega)$, where
$\Omega$ is an arbitrary infinite set, has this property, and it is from this
paper that the term \emph{Bergman's property} arose.

%%%%%%%%%%%%%%%%%%%%%%%%%%%%%

\begin{thm}[cf. \cite{Droste2005aa}]
  \label{thm-scf-cf-bergman}
  Let $G$ be a group. Then $\scf(G)>\aleph_0$ if and only if $\cf(G)>\aleph_0$
  and $G$ has Bergman's property. 
\end{thm}

Droste and G\"obel \cite{Droste2005ab} provide a sufficient criterion for
certain permutation groups to have uncountable strong cofinality and hence Bergman's
property. Their criterion applies to the symmetric group, homeomorphism groups
of the Cantor space, the rationals, and irrationals.

%%%%%%%%%%%%%%%%%%%%%%%%%%%%%

The next lemma connects the notions just defined to that of having a universal
sequence. It appears as Lemma 2.4 in \cite{Maltcev2009aa} and is based on ideas
in Bergman \cite{Bergman2006aa} and Khelif \cite{Khelif2006aa}. 

%%%%%%%%%%%%%%%%%%%%%%%%%%%%%

\begin{lem}
  \label{lem-strong-distortion}
  Let $G$ be a non-finitely generated group and suppose that there exists a
  sequence $\Sq{l_n}$ of natural numbers and an $N\in \N$ such that every
  sequence $\Sq{g_n}$ in $G$ is contained in an $N$-generated subgroup of $G$
  and, for every $n\in \N$, there is a product of length $l_n$ equal to $g_n$ 
  over the $N$ generators. Then $\scf(G)>\aleph_0$.
\end{lem}

Having a universal sequence over a finite alphabet is a stronger property than
the condition in Lemma \ref{lem-strong-distortion}. Hence we obtain the following
corollary. 

%%%%%%%%%%%%%%%%%%%%%%%%%%%%%

\begin{cor}\label{cor-scf}
  If $G$ is a non-finitely generated group and $G$ has an $m$-letter universal
  sequence for some $m\in \N$, then $\scf(G)>\aleph_0$.
\end{cor}

%%%%%%%%%%%%%%%%%%%%%%%%%%%%%

Bergman's original theorem, that the symmetric group has Bergman's property,
follows immediately from Galvin's Theorem \ref{thm-galvin},  Corollary
\ref{cor-scf}, and Theorem \ref{thm-scf-cf-bergman}. It is also possible to
obtain the results of Gourion \cite{Gourion1992aa} and Droste and Holland
\cite{Droste2005aa} as a corollary of our Main Theorem.

%%%%%%%%%%%%%%%%%%%%%%%%%%%%%

\begin{cor}\label{cor-bergman}
  $\scf(\autq)>\aleph_0$ and so $\cf(\autq)>\aleph_0$ and $\autq$ has Bergman's
  property. 
\end{cor}
\begin{proof} 
  This follows immediately from the Main Theorem, Corollary \ref{cor-scf},
  and Theorem \ref{thm-scf-cf-bergman}.
\end{proof}

%%%%%%%%%%%%%%%%%%%%%%%%%%%%%

The ordered set $(\Q, \leq)$ is an example of a \emph{Fra\"iss\'e limit}, namely the limit
of the class of finite linear orders. Automorphism groups of Fra\"iss\'e limits
have many interesting properties; see, for example, \cite{Kechris2007aa},
\cite{Truss1992aa}, and \cite{Truss2003aa}. It is natural to ask if results
analogous to those obtained here for $\autq$ hold for the automorphism groups of
different Fra\"iss\'e limits.

%%%%%%%%%%%%%%%%%%%%%%%%%%%%%

\begin{que} 
  Let $R$ denote the countably infinite random graph. Is it true that every
  countable set of automorphisms of $R$ is contained in an $N$-generated
  subsemigroup or subgroup for some fixed $N\in \N$?  Does $\aut(R)$ have a
  universal sequence over a finite alphabet? Or, more generally, is it possible
  to characterise those Fra\"iss\'e limits whose automorphism groups have either
  of these properties?
\end{que}
%JDM: are there any groups with Sierpinski rank > 2?

%%%%%%%%%%%%%%%%%%%%%%%%%%%%%

The paper is organised as follows: in the next section we provide the relevant
definitions and some general results about universal sequences, and
order-automorphisms of $\Q$. In Section \ref{section-8-letter}, we reduce the
problem of finding a universal sequence for $\autq$ to that of finding a
universal sequence for a subgroup; we also show that $\autq$ has an $8$-letter
universal sequence. We prove our Main Theorem in Section
\ref{section-main-theorem}. 

%%%%%%%%%%%%%%%%%%%%%%%%%%%%%
%%%%%%%%%%%%%%%%%%%%%%%%%%%%%

\section{Preliminaries}

In this section we give the relevant definitions and some results about
universal sequences for arbitrary semigroups. 

In this paper, the natural numbers do not contain $0$, i.e. $\N=\{1,2,3,\dots\}$. For $m\in \N$, let $m\N=\set{mn}{n\in \N}$. The identity function on $\Q$ is denoted by $\id$.
If $f, g\in \autq$, then we define
$$\|f-g\|_{\infty}=\sup\set{|(x)f-(x)g|}{x\in \Q}\in \R\cup\{\infty\}.$$
We denote the \emph{conjugate} $f^{-1}gf$ by $g^f$, and the \emph{commutator}
$f^{-1}g^{-1}fg$ of $f$ and $g$ by $[f,g]$. The \emph{support} of $f\in \aut(\Q,
\leq)$ is
defined and denoted as:
$$\supp(f)=\set{x\in \Q}{(x)f\not=x}$$
and the \emph{fix} of $f$ is just $\fix(f)=\Q\setminus \supp(f)$.  If $X$ is a
subset of $\Q$, then we define the \emph{pointwise stabiliser} of $X$ in $\autq$
by
$$\stab(X)=\set{f\in \aut(\Q, \leq)}{X\subseteq\fix(f)}.$$
The \emph{restriction} of an order-automorphism $f$ to a set $X$ is denoted
$f|_X$. If $f$ setwise stabilises $X$, i.e. $(X)f=X$, then $f|_X\in \aut(X, \leq)$. 
It is well-known that every element of $\autq$ (or any Polish group with a
comeagre conjugacy class) is a commutator; \cite[Theorem 2F]{Glass1981aa},
\cite{Holland1963aa}, \cite{Kechris2007aa}, \cite{Truss1992aa}. 

If $U$ is a set, then we denote the set of sequences of elements of $U$
by $U^\N$. If additionally $U$ is a group, then we use $U^\N$ to denote the
set of sequences of elements of $U$ with componentwise multiplication. 

%%%%%%%%%%%%%%%%%%%%%%%%%%%%%

\begin{prop} \label{prop-general-1}
  Let $S$ be a semigroup, let $U$ and $V$ be subsets of $S$, and let $\Sq{u_n}$,
  $\Sq{v_n}$ be sequences of words over some alphabet $A$.  Then the following hold:
  \begin{enumerate}[\rm (i)]
  \item
    if $\Sq{u_n}$ is universal for $U$ as a subset of $S$, then so is every
    subsequence of $\Sq{u_n}$; 
  
  \item $\Sq{u_n}$ is universal for $U$ as a subset of $S$ if and only if $\Sq{u_n}$ is
     universal for $U^\N$ as a subset of $S^\N$; 

  \item  if $\Sq{u_n}$ is universal for $U$ as a subset of $S$ and
    there exists $B\subseteq S$, and $w_1, w_2, \ldots, w_{m+1}$ in the
    subsemigroup generated by $B$ such that for every $v\in V$, $v=w_1y_1w_2y_2
    \dots w_my_mw_{m+1}$ for some $y_1,y_2, \dots, y_m \in U$, then 
    $$\Sq{w_1u_{m(n-1)+1}w_2\cdots w_mu_{m(n-1)+m}w_{m+1}}$$
    is a universal sequence over $A\cup B$ for $V$ as a subset of $S$;

  \item
    if $m \in \N$, $\Sq{u_n}$ is universal for $U$ as a subset of $S$, and $\Sq{x_n}$ is a sequence
    such that for every $n \in \N$ there are $y_{n,1},y_{n,2}, \dots, y_{n,mn}
    \in U$ where $x_n=y_{n,1}y_{n,2} \cdots y_{n,mn}$, then there is a
    homomorphism $\phi:A^+\to S$ such that
    $$\left(\prod^{mn(n-1)/2+mn}_{i=mn(n-1)/2+1}u_i\right) \phi=x_n$$
    for all $n\in \N$.

\end{enumerate} 
\end{prop} 
\begin{proof}
  {\bf (i).} This is straightforward to verify. \vspace{\baselineskip}

  {\bf (ii).} ($\Rightarrow$) For every $n\in \N$, let
  $\mathbf{x}_n=(x_{m,n})_{m\in\N}$ be a sequence of elements in $U$, i.e.
  $(\mathbf{x}_n)_{n\in \N}$ is a sequence of elements of $U^\N$.  Then, by the
  assumption applied to $(x_{m,n})_{n\in\N}$, for every $m\in \N$ there exists a
  homomorphism $\phi_m:A^{+}\to U$ such that $(u_n)\phi_m=x_{m,n}$. We define
  $\phi:A^{+}\to U^{\N}$ by $$(w)\phi=((w)\phi_1, (w)\phi_2, \ldots).$$ Then
  $\phi$ is a homomorphism and $$(u_n)\phi=(x_{1, n}, x_{2,n},
  \ldots)=\mathbf{x}_n$$ for all $n\in \N$, as required.
  \vspace{\baselineskip}
  
  ($\Leftarrow$) Let $\pi_1:S^{\N}\to S$ be defined by
  $\big((s_n)_{n\in\N}\big)\pi_1=s_1$. Then $\pi_1$ is a homomorphism.  If
  $(x_n)_{n\in\N}$ is a sequence in $U$, then by assumption there exists a
  homomorphism $\phi:A^{+}\to S^{\N}$ such that for all $n$ 
$$(u_n)\phi=(x_n, x_n,
  \ldots).$$ Hence $\phi\pi_1:A^{+}\to S$ is a homomorphism and
  $$(u_n)\phi\pi_1=x_n$$ for all $n\in\N$, as required.\vspace{\baselineskip}
 
  {\bf (iii).} Let $\Sq{x_n}$ be a sequence of elements of $V$. Then for every $n \in \N$ there exist $y_{n,1}, y_{n,2},\dots, y_{n,m} \in U$ such that
$$x_n=w_1y_{n,1}w_2y_{n,2}\dots w_m y_{n,m}w_{m+1}.$$ 
Since $(u_n)_{n\in \N}$ is universal for $U$ as a subset of $S$, there exists a homomorphism
$\phi:A^+\to S$ such that $(u_n)\phi$ is the $n$th element of the sequence $(y_{1,1}, \dots, y_{1,m}, y_{2,1}, \dots, y_{2,m}, \dots)$, i.e.  $(u_{(i-1)m+j})\phi=y_{i,j}$ for all $i,j \in \N$ where $1\leq j \leq m$. The homomorphism $\phi$ can be extended to the natural homomorphism $\Phi:{(A \cup B)}^+ \to S$ satisfying $(b)\Phi=b$ for all $b\in B$. Then, for every $n\in \N$,
\begin{align*}
& (w_1u_{m(n-1)+1}w_2\cdots w_mu_{m(n-1)+m}w_{m+1})\Phi \\
= & (w_1)\Phi(u_{m(n-1)+1})\Phi (w_2)\Phi\cdots (w_m)\Phi (u_{m(n-1)+m})\Phi (w_{m+1})\Phi\\
= & w_1 y_{n,1}w_2\dots w_m y_{n,m}w_{m+1}=x_n,
\end{align*}
as required.
  \vspace{\baselineskip}

  {\bf (iv).} Let $\Sq{x_n}$ be a sequence such that for every $n \in \N$
  there are $y_{n,1}, \ldots, y_{n,mn}\in U$ such that $x_n=y_{n,1} \cdots
  y_{n,mn}$. Since $\Sq{u_n}$ is universal for $U$ as a subset of $S$, there exists a homomorphism $\phi:A^+\to S$ such that $(u_n)\phi$ is the $n$th element of the sequence
$(y_{1,1}, \dots, y_{1,m}, y_{2,1}, \dots, y_{2,2m},y_{3,1}, \dots, y_{3,3m}, \dots)$. 
In other words, $y_{i,j}=(u_{mi(i-1)/2+j})\phi$ for all $i,j\in \N$ where $1\leq j \leq mi$, since $\sum_{k=1}^{i-1}mk=mi(i-1)/2$.
Now
  $$\left(\prod_{k=mn(n-1)/2+1}^{mn(n-1)/2+mn}u_k\right)\phi
  =\prod_{k=mn(n-1)/2+1}^{mn(n-1)/2+mn}(u_k)\phi
=y_{n,1}\cdots
  y_{n,mn}=x_n,$$ as required.  
\end{proof}
 
%%%%%%%%%%%%%%%%%%%%%%%%%%%%%

Note that Proposition \ref{prop-general-1}(ii) holds for arbitrary cartesian
products as well as countable ones. The proof of this result is similar to the
proof in the countable case, but we will not use the more general statement, and
so we have limited ourselves to the countable case.

%%%%%%%%%%%%%%%%%%%%%%

\begin{lem}\label{lem-min-max}
  Let $f,g\in \aut(\Q, \leq)$ be arbitrary, let $\min\{f,g\}$ and $\max\{f,g\}$
  denote the pointwise minimum and maximum of $f$ and $g$, respectively. Then
  $\min\{f,g\}, \max\{f,g\}\in \aut(\Q, \leq)$ and
  ${(\min\{f,g\})}^{-1}=\max\{f^{-1}, g^{-1}\}$.
\end{lem}
\begin{proof}
  Let $h=\min\{f,g\}$, let $k=\max\{f^{-1}, g^{-1}\}$, and let $x,y\in \Q$ be
  such that $x<y$. Then $(x)h=\min\{(x)f, (x)g\}< \min\{(y)f, (y)g\}=(y)h$,
  and so $\min\{f,g\}$ is order-preserving. 

  Suppose without loss of generality that $(x)f\leq (x)g$. Then $(x)fg^{-1}\leq
  x$ and so $(x)hk=(x)fk=\max\{x, (x)fg^{-1}\}=x$. A similar argument shows that
  $(x)kh=x$, and so $h$ and $k$ are bijections, as required.
\end{proof}

%%%%%%%%%%%%%%%%%%%%%%

\section{A reduction of the problem}\label{section-8-letter}

In this section, we prove four lemmas which reduce the problem of finding a
universal sequence for $\autq$ to that of finding such a sequence for a
subgroup.  We also show that the order-automorphisms of the rationals have an
$8$-letter universal sequence.

The first reduction involves the bounded automorphisms, for which we define:
$$B_r=\{g\in \aut(\Q, \leq)\::\: \|g-\id\|_{\infty}\leq r\}$$ 
for all $r\in \Q$, $r>0$.  Note that $B_r$ is not a subgroup of $\autq$ for any
$r\in \Q$, $r>0$. 

%%%%%%%%%%%%%%%%%%%%%%

\begin{lem}\label{lem-step-1}
  Let $(g_{n})_{n\in\N}$ be any sequence in $\aut(\Q, \leq)$. Then there exists
  $p\in \aut(\Q, \leq)$ such that ${g_{n}}^{p}\in B_{2n}$ for all $n\in \N$.
\end{lem} 
\begin{proof}
  We will specify $p$ in terms of the following injective order preserving
  mapping $\sigma:\Z\to \Z$.  We define $\sigma$ recursively, starting with
  $(0)\sigma=0$. Suppose that $\sigma$ is defined on $m\in \Z$. If $m>0$, then
  define $(m+1)\sigma$ to be any value in $\Z$ such that 
  $$(m+1)\sigma>\max\{(m)\sigma g_n,\ (m)\sigma {g_{n}}^{-1}:1\leq n\leq m\};$$ 
  and if $m<0$, then choose $(m-1)\sigma\in \Z$ such that 
  $$(m-1)\sigma <\min\{(m)\sigma g_n,\ (m)\sigma {g_{n}}^{-1}:1\leq n\leq -m\}.$$ 
  
  We will show that no more than one point from $\left(\Z \setminus [-(n-1),n-1]
  \right) \sigma$ lies between $x$ and $(x)g_n$ for all $x\in \Q$ and for all
  $n\in\N$.  Let $x\in \Q$ and $n\in \N$ be arbitrary and suppose that there
  exists $y\in \Z\setminus [-(n-1), n-1]$ such that $(y)\sigma$ lies between $x$
  and $(x)g_n$.  There are four cases to consider depending on the signs of
  $x-(x)g_n$ and $y$.  If $x<(x)g_n$ and $n\leq y$, then $x<(y)\sigma$ and so
  $(x)g_n<(y)\sigma g_n$. By the definition of $\sigma$ and since $n\leq y$,
  $(y+1)\sigma>(y)\sigma g_n$. Hence $(x)g_n<(y)\sigma g_n<(y+1)\sigma$ and so
  $(y+1)\sigma$ is not between $(x)g_n$ and $x$. Thus $\sigma$ has the required
  property. The other cases are analogous.

  Let  $p\in \aut(\Q, \leq)$ be an extension of the function $\sigma:\Z\to \Z$
  and let $x\in \Q$ be arbitrary. There are two cases to consider: $(x)p\leq
  (x)pg_n$ and $(x)p\geq (x)pg_n$. We will only give the proof in the first
  case, the proof of the second case follows by a similar argument. 

  For every $n\in \N$, there exists at most one  $y\in \Z\setminus [-(n-1),
  n-1]$ such that $(y)p=(y)\sigma$ lies between $(x)p$ and $(x)pg_n$. Since
  $p\in \autq$,  it follows that $(x)p\leq (y)p\leq (x)pg_n$ if and only if
  $x\leq y\leq (x)pg_np^{-1}$.  In other words, there is at most one $y\in
  \Z\setminus [-(n-1), n-1]$  that  lies between $x$ and $(x){g_n}^{p^{-1}}$ for
  all $x\in \Q$ and $n\in \N$. Therefore $|(x){g_n}^{p^{-1}}-x| \leq 2n$ for all
  $x\in \Q$ and $n \in \N$.  
\end{proof}

%%%%%%%%%%%%%%%%%%%%%%

\begin{lem}\label{lem-step-2}
  Let $r\in \Q$, $r>0$, and $n\in\N$ be arbitrary. Then $B_{rn}\subseteq
  {(B_r)}^{n}$.
\end{lem}
\begin{proof}
  Let $g\in B_{rn}$ be arbitrary. 
  We use induction on $n$ to show that there exist $h_1 \in B_r$ and
  $h_2 \in B_{r(n-1)}$ with $g={h_1}^{-1}h_2$. We define $h_1$ piecewise as follows:
  \begin{equation*}
    (x)h_1=
    \begin{cases}  
      x-r        & \text{if }  (x)g^{-1}  \leq x-r           \\
      (x)g^{-1}  & \text{if }  x-r     <    (x)g^{-1} < x+r  \\ 
      x+r        & \text{if }  x+r     \leq (x)g^{-1}.
    \end{cases}
  \end{equation*} 
  Then $$h_1=\max\{\min\{x\mapsto x+r, g^{-1}\}, x\mapsto x-r\}$$ and so, by
  Lemma \ref{lem-min-max}, $h_1\in \aut(\Q, \leq)$.  It is clear from the
  definition that $h_1\in B_r$ and so ${h_1}^{-1}\in B_r$. 

  We will show that $h_2=h_1g\in B_{r(n-1)}$. Let $x\in \Q$ be arbitrary. 
  There are three cases to consider. If $(x)g^{-1}\leq x-r$, then
  $x\leq (x-r)g$ and so $(x)h_1g=(x-r)g\geq x$. But $g\in B_{rn}$ and so
  $x-r+rn\geq (x-r)g$. Hence $x \leq (x)h_1 g \leq x + r(n-1)$.  If
  $x-r<(x)g^{-1}<x+r$, then $(x)h_1g=(x)g^{-1}g=x$. In the final case, it
  follows, by a symmetric argument to the first case, that $x-r(n-1)\leq
  (x)h_1g\leq x$ and so $h_1g\in B_{r(n-1)}$.
\end{proof}

%%%%%%%%%%%%%%%%%%%%%%

If $m,n\in \Z$, then we denote the set $\set{mi+n}{i\in \Z}$ by $m\Z+n$ and the
pointwise stabiliser of $m\Z+n$ in $\aut(\Q, \leq)$ by $\stab(m\Z+n)$.

%%%%%%%%%%%%%%%%%%%%%%

\begin{lem}\label{lem-step-3}
  Let $f\in \autq$ be defined by $(x)f=x+1$. Then 
  $$B_{1/3}\subseteq \stab(2\Z)\cdot {\stab(2\Z)}^{f}.$$
\end{lem}
\begin{proof} 
  Let $g\in B_{1/3}$ be arbitrary. Then $2n+2/3\leq (2n+1)g^{-1}\leq 2n+4/3$.
  The closed interval $[2n+2/3, 2n+4/3]$ is a subset of the open interval
  $(2n, 2n+2)$ for all $n\in \Z$. Therefore there is $h\in \stab(2\Z)$ such
  that $(2n+1)h^{-1}=(2n+1)g^{-1}$ for all $n\in \Z$. Since
  ${\stab(2\Z)}^{f}=\stab(2\Z+1)$ and $(2n+1)h^{-1}g=2n+1$ for all $n\in\Z$, it
  follows that $h^{-1}g\in {\stab(2\Z)}^{f}$. Thus $g=h\cdot h^{-1}g\in
  \stab(2\Z)\cdot {\stab(2\Z)}^{f}$, as required.  
\end{proof}

%%%%%%%%%%%%%%%%%%%%%%

If $n\in 2\N$, $n>2$, then we define 
$$I_{n}=\bigcup_{i\in \Z}(ni+2, ni+n)\cap \Q.$$

%%%%%%%%%%%%%%%%%%%%%%

\begin{lem}\label{lem-step-4}
 Let $n\in 2\N$, $n>2$. If $f\in \autq$ is defined by $(x)f=x+1$ , then 
  $$\stab(2\Z)\subseteq \prod_{i=1}^{n/2}{\stab(I_n)}^{f^{2i}}.$$
\end{lem}
\begin{proof}
  Let $h\in \stab(2\Z)$ be arbitrary and for every $i\in \{1, \ldots, n/2\}$
  define $k_i\in \autq$ by
  \begin{equation*}
    (x)k_i=
    \begin{cases}
      (x)h &\text{if }x    \in [nj+2i,nj+2i+2],\ j\in \Z \\
      x    &\text{otherwise.}
    \end{cases}
  \end{equation*}
  Then clearly $h=k_1\cdots k_{n/2}$ and ${k_i}^{f^{-2i}}\in \stab(I_n)$ for all
  $i$, as required. 
\end{proof}

%%%%%%%%%%%%%%%%%%%%%%

In the following corollary we show how the previous four lemmas can be used to
reduce the problem of finding a universal sequence for $\autq$ to that of
finding a universal sequence for $\stab(I_{m})$ as a subset of $\autq$ for any $m\in 2\N$, $m>2$.

%%%%%%%%%%%%%%%%%%%%%%

\begin{cor}\label{cor-combine}
 Let $m\in 2\N$, $m>2$, and let $A$ be an alphabet. If there exists a universal sequence over $A$ for $\stab(I_{m})$ as a subset of $\autq$, then there is a universal sequence for $\autq$ over
  $A\cup \{f, f^{-1}\}$. 
\end{cor}
\begin{proof} 
  It follows from Proposition \ref{prop-general-1}(iii) and Lemmas 
  \ref{lem-step-3} and \ref{lem-step-4} that if there is a
  universal sequence over $A$ for $\stab(I_{m})$ as a subset of $\autq$, then there is a
  universal sequence $\Sq{u_n}$ for $B_{1/3}$ as a subset of $\autq$ over the alphabet $A\cup \{f, f^{-1}\}$.
We define the sequence $\Sq{w_n}$ over $A\cup \{f, f^{-1}\}$ by
$$w_n=\prod_{i=6n(n-1)/2+1}^{6n(n-1)/2+6n}u_i$$ 
for all $n\in \N$. 

 Suppose that $\Sq{g_n}$ is an arbitrary sequence in $\autq$. Then,
  by Lemma \ref{lem-step-1}, there exists $p\in \autq$ such that ${g_n}^p\in
  B_{2n}$ for all $n\in\N$.
Since $B_{2n}\subseteq {(B_{1/3})}^{6n}$ for all $n\in\N$ (by Lemma \ref{lem-step-2}), it follows from Proposition \ref{prop-general-1}(iv) that there exists a 
  homomorphism $$\phi:{\big(A\cup \{f, f^{-1}\}\big)}^+\to \autq$$ such that
  $$(w_n)\phi={g_n}^p$$
  for all $n\in \N$.  Conjugating by
  $p^{-1}$ is an automorphism of $\autq$ and so
  $$\theta:{\big(A\cup \{f, f^{-1}\}\big)}^{+}\to \autq$$
  defined by $(x)\theta={((x)\phi)}^{p^{-1}}$ is a homomorphism and
  $(w_n)\theta=g_n$ for all $n\in \N$. In other words, $\Sq{w_n}$ is a universal
  sequence for $\autq$ over $A\cup \{f, f^{-1}\}$, as required.
\end{proof}

%%%%%%%%%%%%%%%%%%%%%%

We will prove the Main Theorem in the next section.
Since it is significantly more complicated to prove that there is a $2$-letter
universal sequence for $\aut(\Q, \leq)$, we first show that there is such a
sequence over an $8$-letter alphabet. Thus for any reader who is only interested
to learn that there is a universal sequence over a finite alphabet, the next
theorem ought to suffice.

%%%%%%%%%%%%%%%%%%%%%%

\begin{thm} \label{thm-8}
  There is a universal sequence over an $8$-letter alphabet for $\aut(\Q,
  \leq)$.  
\end{thm}
\begin{proof}
  By Corollary \ref{cor-combine}, it suffices to find a $6$-letter universal
  sequence for $\stab(I_{4})$ as a subset of $\autq$.  We set $\Omega=\Z\times\Z\times\Q$ equipped with
  the usual lexicographic order. Then since $\Omega$ is a countable dense linear
  order without endpoints, it follows that $\Omega$ is order-isomorphic to $\Q$.
  We consider $\aut(\Omega, \leq)$ instead of $\autq$ in this proof.  Set
  $\Omega_0=\{0\}\times\{0\}\times \Q$ and 
  $$\stab(\Omega\setminus \Omega_0):=\set{g\in \aut(\Omega, \leq)}{\supp(g)\subseteq
  \Omega_0}.$$
  Since $\Omega_0$ is a bounded open interval in $\Omega$, there is an
  order-isomorphism $$\theta:(4i-1, 4i+3)\cap \Q\to \Omega$$ for all $i\in \Z$
  such that $$\left((4i, 4i+2)\cap \Q\right)\theta=\Omega_0.$$  Conjugation by
  $\theta$ is a group-isomorphism from $\aut((4i-1, 4i+3)\cap \Q, \leq)$ to
  $\aut(\Omega, \leq)$ mapping $\stab(((4i-1,4i+3)\setminus (4i, 4i+2)))\cap \Q)$ to
  $\stab(\Omega\setminus \Omega_0)$. It follows that there is a group-isomorphism from $\stab(4\Z-1)$ to the direct product ${\aut(\Omega,
  \leq)}^{\Z}$ mapping $\stab(I_{4})$ to ${\stab(\Omega\setminus \Omega_0)}^{\Z}$.
  Therefore it suffices, by Proposition \ref{prop-general-1}(ii), to show that
  there is a $6$-letter universal sequence for $\stab(\Omega\setminus \Omega_0)$
  as a subset of $\aut(\Omega,\leq)$.

  Let $\Sq{g_n}$ be such that $g_n\in \stab(\Omega\setminus \Omega_0)$ for all
  $n\in\N$. Then, since every element of $\aut(\Omega, \leq)$ is a commutator,
  there exist $h_{2n-1}, h_{2n}\in\stab(\Omega\setminus \Omega_0)$ such that
  $g_n=[h_{2n-1}, h_{2n}]$ for all $n\in\N$. Since $\supp\left(h_{n}\right)\subseteq
  \Omega_0$, for every $n\in \N$, we can define $\bar{h}_{n}\in \autq$ such that
  $(0,0,x)h_n=(0, 0, (x)\bar{h}_{n})$ for all $x\in \Q$. 
  
  We define $a,b,c \in \aut(\Omega, \leq)$ by
  \begin{equation*}
    (i,j,x)a=
    \begin{cases} 
      (i,j,(x)\bar{h}_{2n})   &\text{if }i=-2n,\ j=0\\
      (i,j,(x)\bar{h}_{2n-1}) &\text{if }i=2n-1,\ j=0\\
      (i,j,x)                 &\text{otherwise,}
    \end{cases}
  \end{equation*}
  \begin{equation*}
    (i,j,x)b=(i+1,j,x)
    \qquad\text{and}\qquad
    (i,j,x)c=
    \begin{cases}
      (i,j+1,x) &\text{if }i\not=0\\
      (i,j,x)   &\text{if }i=0
    \end{cases}
  \end{equation*}
  for all $(i,j,x)\in \Omega$. It is routine to verify that
  $a,b,c\in\aut(\Omega, \leq)$.

  If $(0, 0, x)\in \Omega_0$ is arbitrary, then 
  \begin{eqnarray*}
    (0,0,x)a^{b^{1-2n}}&=&(2n-1,0,x)ab^{1-2n}=(2n-1,0, (x)\bar{h}_{2n-1})b^{1-2n}=
    (0,0, (x)\bar{h}_{2n-1})\\
    &=&(0,0,x)h_{2n-1}\in \Omega_0
  \end{eqnarray*}
  and 
  $$(0,0,x)a^{b^{2n}c}=(-2n,0,x)ab^{2n}c=(-2n,0,(x)\bar{h}_{2n})b^{2n}c=
  (0,0, (x)\bar{h}_{2n})=(0,0,x)h_{2n}\in \Omega_0.$$
  Hence on $\Omega_0$, at least, $[a^{b^{1-2n}}, a^{b^{2n}c}]$ equals
  $[h_{2n-1}, h_{2n}]=g_n$.

  Since $\supp(a)\subseteq \Z\times\{0\}\times \Q$,  
  $$\supp\left(a^{b^{1-2n}}\right)=\supp(a)b^{1-2n}\subseteq \Z\times \{0\}\times \Q$$
  and
  $$\supp\left(a^{b^{2n}c}\right)\subseteq (\Z\times\{1\}\times\Q)\cup \Omega_0.$$
  Thus $\supp\left(a^{b^{1-2n}}\right)\cap \supp\left(a^{b^{2n}c}\right)\subseteq
  \Omega_0$. Since $a^{b^{1-2n}}$ and $a^{b^{2n}c}$ also fix $\Omega_0$ setwise it follows that 
$\supp\left(\left[a^{b^{1-2n}}, a^{b^{2n}c}\right]\right)\subseteq \Omega_0$.  
  Hence, for all $n\in \N$, 
$$\left[a^{b^{1-2n}}, a^{b^{2n}c}\right]=\left[h_{2n-1}, h_{2n}\right]=g_n.$$ 
The map that takes each letter in the alphabet $A=\{a, a^{-1}, b, b^{-1}, c, c^{-1}\}$ to the corresponding element of $\autq$ defined above extends to a unique homomorphism $\phi:A^{+}\to \autq$ and we have just shown that $([a^{b^{1-2n}}, a^{b^{2n}c}])\phi=g_n$ for all $n \in \N$.

Thus, since $\Sq{g_n}$ was an arbitrary sequence of elements of $\stab(\Omega\setminus \Omega_0)$, it follows that the sequence
\begin{equation}
\label{IM_universal}
\Sq{[a^{b^{1-2n}}, a^{b^{2n}c}]}
\end{equation}
is universal for  $\stab(\Omega\setminus \Omega_0)$ as a subset of $\autq$ over the $6$-letter alphabet $A$, which concludes the proof.
\end{proof}

%%%%%%%%%%%%%%%%%%%%%%

The proof of Theorem \ref{thm-8} establishes the existence of an 8-letter
universal sequence for $\autq$; we will now construct such a sequence
explicitly.

Let $m \in 2\N$ with $m>2$ be fixed, let $\Sq{w_n}$ be a universal sequence for $\stab(I_m)$ as a subset of $\autq$ over some alphabet $A$, and let $\Sq{g_n}$ be an arbitrary sequence of elements of $\autq$.
By Lemma \ref{lem-step-1} there exists $p\in \autq$ such that ${g_n}^p \in B_{2n}$ for all $n\in \N$. Hence, by Lemma \ref{lem-step-2}, there exist $u_{n,1}, u_{n,2}, \dots, u_{n,6n}\in B_{1/3}$ such that
\begin{equation}\label{8combine1}
{g_n}^p=\prod_{i=1}^{6n} u_{n,i}
\end{equation}
for all $n\in \N$. Let $f\in \autq$ be defined by $(x)f=x+1$.  Then, by Lemma \ref{lem-step-3}, there exist  $v_{n, i,1}, v_{n,i,2}$ in $\stab(2\Z)$ auch that 
\begin{equation}\label{8combine2}
u_{n,i}= v_{n, i,1} \cdot {v_{n,i,2}}^{f}
\end{equation}
for all $n\in \N$ and $1\leq i\leq 6n$. By Lemma \ref{lem-step-4}, there exist 
$t_{n,i,j,1}, t_{n,i,j,2}\dots, t_{n,i,j,m/2} \in \stab(I_m)$ such that
\begin{equation}\label{8combine3}
v_{n,i,j}=\prod_{k=1}^{m/2}{t_{n,i,j,k}}^{f^{2k}}
\end{equation}
for all $n\in \N$, $1\leq i\leq 6n$, $1\leq j\leq 2$.

Finally, combining equations \eqref{8combine1}, \eqref{8combine2} and \eqref{8combine3} above we have that
\begin{equation}\label{8combine4}
{g_n}^p=\prod_{i=1}^{6n} u_{n,i}=\prod_{i=1}^{6n}\left( v_{n, i,1} \cdot {v_{n,i,2}}^{f}\right)
=\prod_{i=1}^{6n}\left( \prod_{k=1}^{m/2} {t_{n,i,1,k}}^{f^{2k}}\cdot \prod_{k=1}^{m/2}{t_{n,i,2,k}}^{f^{2k+1}}\right)
\end{equation}
for all $n\in \N$.

The set 
$$T=\set{t_{n,i,j,k}}{n\in \N, 1\leq i\leq 6n, 1\leq j\leq 2, 1\leq k\leq m/2}$$ 
is contained in $\stab(I_m)$ and $\Sq{w_n}$ is universal for $\stab(I_m)$ as a subset of $\autq$. Hence, if the elements of $T$ are ordered in any way, then there exists a homomorphism $\phi:A^+\to \autq$ that maps $w_n$ to the $n$th element of $T$. Specifically, if we order $T$ according to the usual lexicographical order on the tuples $(n,i,j,k)$, then $t_{n,i,j,k}$ is in position 
$$3mn(n-1)+(i-1)m+\frac{(j-1)m}{2}+k=(n,i,j,k)\iota.$$
In other words,
\begin{equation*}
\left(w_{(n,i,j,k)\iota}\right)=\left(w_{3mn(n-1)+(i-1)m+\frac{(j-1)m}{2}+k}\right)\phi=t_{n,i,j,k}.
\end{equation*}
We may extend $\phi$ to a homomorphism $\phi_1: {(A\cup \{f, f^{-1}\})}^+\to \autq$ mapping the letters $f$ and $f^{-1}$ to the automorphisms $f$ and $f^{-1}$, respectively. Then by \eqref{8combine4} above it follows that
\begin{align*}
&\left( \prod_{i=1}^{6n}\left( \prod_{k=1}^{m/2} {w_{(n,i,1,k)\iota}}\cdot \prod_{k=1}^{m/2}{w_{(n,i,2,k)\iota}}^{f^{2k+1}}\right)\right)\phi_1\\
= &\prod_{i=1}^{6n}\left( \prod_{k=1}^{m/2} \left({w_{(n,i,1,k)\iota}}^{f^{2k}}\right)\phi_1\cdot \prod_{k=1}^{m/2}\left({w_{(n,i,2,k)\iota}}^{f^{2k+1}}\right)\phi_1\right)\\
= &\prod_{i=1}^{6n}\left( \prod_{k=1}^{m/2} {t_{n,i,1,k}}^{f^{2k}}\cdot \prod_{k=1}^{m/2}{t_{n,i,2,k}}^{f^{2k+1}}\right)
={g_n}^p
\end{align*}
for all $n\in \N$. Conjugation by $p^{-1}$ is an automorphism of $\autq$. Composing $\phi_1$ with this automorphism gives another homomorphism $\phi_2: {(A\cup \{f, f^{-1}\})}^+\to \autq$ and 

$$\left(\prod_{i=1}^{6n}\left( \prod_{k=1}^{m/2} {w_{(n,i,1,k)\iota}}^{f^{2k}}\cdot \prod_{k=1}^{m/2}{w_{(n,i,2,k)\iota}}^{f^{2k+1}}\right)\right)\phi_2={({g_n}^p)}^{p^{-1}}=g_n$$
for all $n\in \N$. Since $\Sq{g_n}$ was an arbitrary sequence of elements of $\autq$, it follows that the sequence with $n$th term equal to
\begin{align*}
&\prod_{i=1}^{6n}\left( \prod_{k=1}^{m/2} {w_{(n,i,1,k)\iota}}^{f^{2k}}\cdot \prod_{k=1}^{m/2}{w_{(n,i,2,k)\iota}}^{f^{2k+1}}\right)\\
=&\prod_{i=1}^{6n}\left( \prod_{k=1}^{m/2} {w_{3mn(n-1)+(i-1)m+k}}^{f^{2k}}\cdot \prod_{k=1}^{m/2}{w_{3mn(n-1)+(i-1)m+\frac{m}{2}+k}}^{f^{2k+1}}\right)\\
=&\prod_{i=3n(n-1)}^{3n(n+1)-1}\left( 
\prod_{k=1}^{m/2}{w_{mi+k}}^{f^{2k}}\cdot 
\prod_{k=1}^{m/2}{w_{\frac{m(2i+1)}{2}+k}}^{f^{2k+1}}
\right)
\end{align*}
 is universal for $\autq$ over the alphabet $A\cup \{f, f^{-1}\}$.

In the case that $m=4$, the $n$th term of this universal sequence is
\begin{align}
 \prod_{i=3n(n-1)}^{3n(n+1)-1}
{w_{4i+1}}^{f^{2}}\cdot {w_{4i+2}}^{f^{4}}\cdot
{w_{4i+3}}^{f^{3}}\cdot {w_{4i+4}}^{f^{5}}.\label{eq-univ-sequence}
\end{align}
Letting $\Sq{w_n}$ be the sequence given by \eqref{IM_universal} now gives the universal sequence
\begin{equation*}
 \Sq{\prod_{i=3n(n-1)}^{3n(n+1)-1}
  [a^{b^{-8i-1}}, a^{b^{8i+2}c}]^{f^{2}}
  \cdot [a^{b^{-8i-3}}, a^{b^{8i+4}c}]^{f^{4}}
  \cdot [a^{b^{-8i-5}}, a^{b^{8i+6}c}]^{f^{3}}
  \cdot [a^{b^{-8i-7}}, a^{b^{8i+8}c}]^{f^{5}}}
\end{equation*}
 for $\autq$  over the alphabet $\{a,a^{-1}, b, b^{-1}, c, c^{-1}, f, f^{-1}\}$.

%%%%%%%%%%%%%%%%%%%%%%
%%%%%%%%%%%%%%%%%%%%%%

\section{Proof of the Main Theorem}\label{section-main-theorem}

In this section we prove the Main Theorem, which we restate for the convenience
of the reader. \medskip

\noindent\textbf{Main Theorem.}
\textit{There is a $2$-letter universal sequence for $\autq$.}\medskip

If $X$ is a totally ordered set, then we denote by $X^*$ the set $X\cup
\{\infty\}$ where the order of $X$ is extended by adjoining a 
maximum element $\infty\not\in X$. 

We identify $\Q$ with the set $\Omega=\Z\times\Z^*\times\Z^*\times\Q^*$ equipped with
the usual lexicographic order. Then, as in the proof of Theorem \ref{thm-8},
$\Omega$ is order-isomorphic to $\Q$. It is straightforward to verify that there
is an order-isomorphism $\phi$ from $\Q$ to $\Omega$ such that 
$\Q\cap (4n-1, 4n+1)$ 
is mapped to $$\Omega_n=\{n\}\times\{0\}\times \{0\}\times (-1,1)$$ and 
$\Q\cap [4n-2, 4n+2]$ is mapped to 
$$\set{\alpha\in\Omega}{(n-1,\infty, \infty, \infty)\leq \alpha\leq(n,\infty,
\infty, \infty)},$$
for every $n\in \Z$.  Moreover, $\phi$ can be chosen so that the function $f$
obtained by conjugating $x\mapsto x+1$ by $\phi$ satisfies
$$(i,j,k,x)f^4=(i+1,j,k,x)$$
for all $(i,j,k,x)\in \Omega$. We will only make use of powers of $f^4$ in the
remainder of the paper, and so we do not require (or give) an explicit
description of the action of $f$ itself on $\Omega$. 

We consider $\aut(\Omega, \leq)$ rather than $\autq$ for the remainder of this
section. 

Let $\Sq{h_n}$ be an
arbitrary sequence of elements in
$\stab(\Omega\setminus\bigcup_{n\in\Z}\Omega_{12n})$.  We will show that there
exists
$g\in \auto$ such that $$\left[{\left(g\cdot g^{f^{-12}}\right)}^{{\left(g^{f^{-4}}\right)}^ng^{f^{-28}}},
{\left(g\cdot g^{f^{-12}}\right)}^{{\left(g^{f^{-4}}\right)}^{-n}}\right]=h_n$$ for all $n\in \N$.  
In other words, $\stab(\Omega\setminus\bigcup_{n\in\Z}\Omega_{12n})$, as a subset of $\aut(\Omega,\leq)$, has a universal
sequence over $\{f,f^{-1},g,g^{-1}\}$. 

Note that, by definition, $f^{\phi^{-1}}=\phi f \phi^{-1}$ is the map $x \mapsto x+1$ and so the order-isomorphism 
$\phi f^{-1}=(\phi f^{-1} \phi^{-1})\phi=(\phi f \phi^{-1})^{-1} \phi$ maps $[48n,48n+2]$ to $\Omega_{12n}$, for all $n \in \Z$. Hence
\begin{align*}
{\stab(I_{48})}^{\phi f^{-1}}
=&f \phi^{-1} \stab(I_{48}) \phi f^{-1}=\stab(\Omega\setminus\bigcup_{n\in\Z}\Omega_{12n}).
\end{align*}
Thus, it will follow that $\stab(I_{48})$, as a subset of $\autq$, has a universal sequence over 
$\{f,f^{-1},g,g^{-1}\}$ and so, by Corollary \ref{cor-combine}, $\autq$ has
a universal sequence over the same alphabet. Once we have defined $g$, we
show in Lemma \ref{lem-group-semigroup} that the group generated by $f$ and $g$
equals the semigroup generated by $f^{-48}g$ and $f$. More precisely, each of
$f$, $g$, $f^{-1}$, and $g^{-1}$ is equal to an explicit product over
$f^{-48}g$ and $f$, which is independent of the sequence $\Sq{h_n}$. Therefore
we will have shown that  there is a universal sequence for $\autq$ over the alphabet
$\{f^{-48}g, f\}$.

%%%%%%%%%%%%%%%%%%%%%%

Since $(h_m)_{m \in \N}$ is a sequence in
$\stab(\Omega\setminus\bigcup_{n\in 12\Z}\Omega_{n})$, it follows that
$h_m|_{\Omega_{n}}\in \aut(\Omega_{n}, \leq\nobreak)$ for all $m \in \N$ and $n\in 12\Z$. 
Every element of $\aut(\Omega_{n}, \leq)\cong\autq$ is a commutator, and 
so there exist $k_{m,n}, k_{-m, n}\in \aut(\Omega_{n}, \leq\nobreak)$ such
that 
$$[k_{-m,n}, k_{m,n}]=h_{m}|_{\Omega_{n}}$$
for all $m\in \N$, $n\in 12\Z$ and we define $k_{0,n}$ to be the identity for
all $n\in 12\Z$. For every $m\in\Z$ and $n\in 12\Z$ there exists
$\bar{k}_{m,n}\in \aut((-1,1)\cap \Q, \leq)$ such that 
$$(n,0,0,y)k_{m,n}=(n,0,0,(y)\bar{k}_{m,n})$$
where $y\in (-1,1)\cap \Q$.

To define the required $g$, we specify four auxiliary order-automorphisms
$a,b,c,d$ of $\Omega$:
\begin{eqnarray*}
  (i,j,m,x)a&=&
  \begin{cases}
    (i,j,m,(x)\bar{k}_{m,i})      &\text{if }i\in 24\Z,\ j\in 2\Z,  \ x\in (-1,1)\\
    (i,j,m,(x){\bar{k}_{m,i}}^{-1}) &\text{if }i\in 24\Z,\ j\in 2\Z+1,\ x\in (-1,1)\\
    (i,j,m,(x)\bar{k}_{-m,i})     &\text{if }i\in 24\Z+12,\ j\in 2\Z,  \ x\in (-1,1)\\
    (i,j,m,(x){\bar{k}_{-m,i}}^{-1})&\text{if }i\in 24\Z+12,\ j\in 2\Z+1,\ x\in (-1,1)\\
    (i,j,m,x)                     &\text{otherwise,}      
  \end{cases}
\end{eqnarray*}
\begin{eqnarray*}
  (i,j,m,x)b&=&
  \begin{cases}
    (i,j+1,m,x) &\text{if }i\in 24\Z\\
    (i,j-1,m,x) &\text{if }i\in 24\Z+12\\
    (i,j,m,x)   &\text{otherwise},
  \end{cases}
  \end{eqnarray*}
\begin{eqnarray*}
  (i,j,m,x)c&=&
  \begin{cases}
    (i,j,m+1,x) &\text{if }i\in 24\Z\\
    (i,j,m-1,x) &\text{if }i\in 24\Z+12\\
    (i,j,m,x)   &\text{otherwise},
  \end{cases}
  \end{eqnarray*}
\begin{eqnarray*}
  (i,j,m,x)d&=&
  \begin{cases}
    (i,j,m,x+2) &\text{if }i\in 24\Z,\    (j,m)\not=(0,0)\\
    (i,j,m,x-2) &\text{if }i\in 24\Z+12,\ (j,m)\not=(0,0)\\
    (i,j,m,x)   &\text{otherwise}.
  \end{cases}
\end{eqnarray*}
It is routine to verify that $a,b,c,d\in \auto$ and that 
$$b^{f^{48}}=b^{-1},\quad c^{f^{48}}=c^{-1},\quad d^{f^{48}}=d^{-1},\quad bc=cb.$$

%%%%%%%%%%%%%%%%%%%%%%

We are now able to define the second automorphism $g\in \autq$ required to
generate $\Sq{h_n}$: 
\begin{equation}\label{eq-functions-inverses}
  g=ab\cdot c^{f^4}\cdot {(b^{-1})}^{f^{12}}\cdot {d}^{f^{28}}.
\end{equation}
Since 
\begin{equation*}
  \begin{array}{ll}
    \supp(ab)\subseteq 12\Z\times \Z^*\times\Z^*\times\Q^*,&
    \supp\left(c^{f^4}\right)\subseteq (12\Z+1)\times \Z^*\times\Z^*\times\Q^*\\
    \supp\left({(b^{-1})}^{f^{12}}\right)\subseteq(12\Z+3)\times\Z^*\times\Z^*\times\Q^*,&
    \supp\left(b^{f^{28}}\right)\subseteq (12\Z+7)\times \Z^*\times\Z^*\times\Q^*,
  \end{array}
\end{equation*}
the supports of $ab$, $c^{f^4}$, ${(b^{-1})}^{f^{12}}$, and $d^{f^{28}}$ are disjoint. In particular, this implies that these automorphisms commute.

%%%%%%%%%%%%%%%%%%%%%%

\begin{lem}\label{abab=bb}
  $(ab)^2=b^2$. 
\end{lem}
\begin{proof}
  From the definitions of $a$ and $b$, we have that:
  \begin{eqnarray*}
    (i,j,m,x)ab&=&
    \begin{cases}
    (i,j+1,m,(x)\bar{k}_{m,i})      &\text{if }i\in 24\Z,\ j\in 2\Z,  \ x\in (-1,1)\\
    (i,j+1,m,(x)\bar{k}_{m,i}^{-1}) &\text{if }i\in 24\Z,\ j\in 2\Z+1,\ x\in (-1,1)\\
    (i,j-1,m,(x)\bar{k}_{-m,i})     &\text{if }i\in 24\Z+12,\ j\in 2\Z,  \ x\in (-1,1)\\
    (i,j-1,m,(x)\bar{k}_{-m,i}^{-1})&\text{if }i\in 24\Z+12,\ j\in 2\Z+1,\ x\in (-1,1)\\
    (i,j,m,x)b                     &\text{otherwise.}      
    \end{cases}
  \end{eqnarray*}
  If $i,j,m,x$ do not fulfil any of the first 4 conditions in the displayed
  equation above, then clearly $(i,j,m,x){(ab)}^2=(i,j,m,x)b^2$.
  If $i\in 24\Z$, $j\in 2\Z$, and $x\in (-1,1)$, then 
  $$(i,j,m,x){(ab)}^2=(i,j+1,m,(x)\bar{k}_{m,i})ab
  =(i,j+2,m,(x)\bar{k}_{m,i}{\bar{k}_{m,i}}^{-1})=(i,j,m,x)b^2.$$
  The remaining cases follows by similar arguments.
\end{proof}

%%%%%%%%%%%%%%%%%%%%%%

\begin{lem}\label{lem-group-semigroup}
  The semigroup generated by $f^{-48}g$ and $f$ is the group generated by $f$
  and $g$.   More precisely, each of $f$, $g$, $f^{-1}$, and $g^{-1}$ is equal
  to a fixed product over $f^{-48}g$ and $f$ which is independent of the
  sequence $\Sq{h_n}$.
\end{lem}
\begin{proof} 
Let $S$ be the semigroup generated by $f^{-48}g$ and $f$. Clearly $f$ and $g=f^{48}(f^{-48}g)$ are in $S$.
We will now show that
${(g^2)}^{f^{48}}=g^{-2}.$

Since $ab$, $c^{f^4}$, ${(b^{-1})}^{f^{12}}$, and
  $d^{f^{28}}$ commute and ${(ab)}^2=b^2$ (Lemma \ref{abab=bb}),
  $$g^2={(ab)}^2{(c^2)}^{f^4}{(b^{-2})}^{f^{12}}{(d^2)}^{f^{28}}
  =b^2{(c^2)}^{f^4}{(b^{-2})}^{f^{12}}{(d^2)}^{f^{28}}.$$ Therefore, by 
  equation \eqref{eq-functions-inverses},
  $${(g^2)}^{f^{48}}={(b^2)}^{f^{48}}{\left({(c^2)}^{f^{48}}\right)}^{f^4}
  {\left({(b^{-2})}^{f^{48}}\right)}^{f^{12}}{\left({(d^2)}^{f^{48}}\right)}^{f^{28}}
  =b^{-2}{(c^{-2})}^{f^4}{(b^{2})}^{f^{12}}{(d^{-2})}^{f^{28}}=g^{-2}.$$ 
Thus 
\begin{equation}\label{g_inverse}
g^{-1}={(g^2)}^{f^{48}}g=(f^{-48}g)f^{48}(f^{-48}g)f^{96}(f^{-48}g)\in S,
\end{equation}
and so 
\begin{equation}\label{f_inverse}
f^{-1}=f^{47}(f^{-48}g)g^{-1}=f^{47}(f^{-48}g)(f^{-48}g)f^{48}(f^{-48}g)f^{96}(f^{-48}g)\in S,
\end{equation}
which concludes the proof.
\end{proof}

%%%%%%%%%%%%%%%%%%%%%%

In the following three lemmas, we show that 
$$\left[{\left(g\cdot g^{f^{-12}}\right)}^{{\left(g^{f^{-4}}\right)}^mg^{f^{-28}}},
{\left(g\cdot g^{f^{-12}}\right)}^{{\left(g^{f^{-4}}\right)}^{-m}}\right]=\left[a^{c^md}, a^{c^{-m}}\right]=h_m,$$
for all $m\in \N$.

Suppose that $u,v\in \auto$ are such that 
$\supp(u), \supp(v)\subseteq 12\Z\times \Z^*\times\Z^*\times\Q^*.$
Then $\supp\left(u^{f^{4i}}\right)\cap \supp\left(v^{f^{4j}}\right)=\emptyset$
for all $i,j\in \Z$ such that $4i\not=4j\pmod{48}$. 
It follows that $u^{f^{4i}}$ and $v^{f^{4j}}$ commute for any such $i,j$, and, in particular, 
this holds when $u$ or $v$ is any product of $a$, $b$,
$c$, or $d$.

%%%%%%%%%%%%%%%%%%%%%%

\begin{lem}\label{step3}
 $\left[{\left(g\cdot g^{f^{-12}}\right)}^{{\left(g^{f^{-4}}\right)}^mg^{f^{-28}}},
{\left(g\cdot g^{f^{-12}}\right)}^{{\left(g^{f^{-4}}\right)}^{-m}}\right]=\left[a^{c^md}, a^{c^{-m}}\right]$
for all
  $m\in \N$.
\end{lem}
\begin{proof} 
  Since $g=abc^{f^4}{\left(b^{-1}\right)}^{f^{12}}d^{f^{28}}$, it follows that 
  $$g^{f^{-12}}={(ab)}^{f^{-12}}c^{f^{-8}}(b^{-1})d^{f^{16}}$$
  and so 
  $$g\cdot g^{f^{-12}}=
  \big(abc^{f^4}{(b^{-1})}^{f^{12}}d^{f^{28}}\big)\big({(ab)}^{f^{-12}}c^{f^{-8}}
  (b^{-1})d^{f^{16}}\big).$$
  Since the only pair $(4i, 4j)$ of powers of $f$ in $\{0,4,12,28\}\times \{-12,
  -8, 0, 16\}$ in this product such that $4i=4j\pmod{48}$ is $(0,0)$, it follows
  that 
  $$g\cdot g^{f^{-12}}={(ab)}^{f^{-12}}c^{f^{-8}}a
  c^{f^4}{\left(b^{-1}\right)}^{f^{12}}d^{f^{16}}d^{f^{28}}.$$
  Also
  $${\left(g^{f^{-4}}\right)}^m={\left({(ab)}^m\right)}^{f^{-4}}c^m {\left(b^{-m}\right)}^{f^8}{\left(d^m\right)}^{f^{24}}.$$
  The only $(4i,4j)$ in $\{-12, -8, 0, 4, 12, 16\}\times \{-4,0,8,24\}$
  such that $4i=4j\pmod{48}$ is $(0,0)$, which implies that
    $${\left(g\cdot g^{f^{-12}}\right)}^{{\left(g^{f^{-4}}\right)}^{-m}}=
    {(ab)}^{f^{-12}}c^{f^{-8}}a^{c^{-m}}c^{f^4}{(b^{-1})}^{f^{12}}d^{f^{16}}d^{f^{28}}.$$
  Next 
  $$g^{f^{-28}}={(ab)}^{f^{-28}}c^{f^{-24}}{(b^{-1})}^{f^{-16}}d$$
  and the only $(4i,4j)$ in $\{-12, -8, 0, 4, 12, 16, 28\}\times \{-28,-24,-16, 0\}$
  such that $4i=4j\pmod{48}$ is $(0,0)$. Therefore
  $${\left(g\cdot g^{f^{-12}}\right)}^{{\left(g^{f^{-4}}\right)}^mg^{f^{-28}}}
  ={(ab)}^{f^{-12}}c^{f^{-8}}a^{c^md}c^{f^4}{(b^{-1})}^{f^{12}}d^{f^{16}}d^{f^{28}}$$
  and hence
 $$\left[{\left(g\cdot g^{f^{-12}}\right)}^{{\left(g^{f^{-4}}\right)}^mg^{f^{-28}}},
{\left(g\cdot g^{f^{-12}}\right)}^{{\left(g^{f^{-4}}\right)}^{-m}}\right]=\left[a^{c^md}, a^{c^{-m}}\right],$$
  as required.
\end{proof}

%%%%%%%%%%%%%%%%%%%%%%

\begin{lem}\label{step1}
  $\left[a^{c^md}, a^{c^{-m}}\right]\big|_{\Omega_{n}}=\left[k_{-m, n}, k_{m,n}\right]|_{\Omega_{n}}
  =h_{m}|_{\Omega_{n}}$ for all $n\in 12\Z$ and $m \in \N$.
\end{lem}
\begin{proof} 
  Let $x\in (-1,1)$ and let $m\in \Z$ be arbitrary. 
  If $n\in 24\Z$, then 
  $(n, 0, 0, x)c^m=(n, 0, m, x),$
  whereas if $n\in 24\Z+12$, then 
  $(n, 0, 0, x)c^m=(n, 0, -m, x).$
  Hence, in either case (i.e. if $n\in 12\Z$) 
  $$(n, 0, 0, x)a^{c^{-m}}=(n, 0, 0, (x)\bar{k}_{m,n})=(n,0,0,x)k_{m,n}$$
  and, similarly, 
  $$(n, 0, 0, x)a^{c^{m}}=(n, 0, 0, x)k_{-m,n}.$$
  Since $d$ fixes the points in $\Omega$ with second and third component equal
  to $0$, it follows that 
  $$(n, 0, 0, x)a^{c^{m}d}=(n, 0, 0, x)k_{-m,n}.$$
  In other words, since $\Omega_{n}=\{n\}\times \{0\}\times \{0\}\times \Q$, it
  follows that 
  $a^{c^{-m}}|_{\Omega_n}=k_{m,n}$ and $a^{c^{m}d}|_{\Omega_n}=k_{-m, n}\in
  \aut(\Omega_n, \leq)$, $n\in 12\Z$.
  Thus $$\left[a^{c^md}, a^{c^{-m}}\right]\big|_{\Omega_n}=\left[a^{c^md}|_{\Omega_n},
  a^{c^{-m}}|_{\Omega_n}\right]=\left[k_{-m, n}, k_{m,n}\right]|_{\Omega_n}=h_{m}|_{\Omega_n},$$
  as required.
\end{proof}

%%%%%%%%%%%%%%%%%%%%%%

We will use the following observation in the proof of the next lemma. If $f$ and $g$ are permutations of a set $X$ and $Y\subseteq X$, then
\begin{equation}\label{conjugate_support}
\supp\left(f^{g}\right)\cap Y=\supp(f)g \cap Y=\left(\supp(f)\cap Yg^{-1}\right)g.
\end{equation}

\begin{lem}
  $\left[a^{c^md}, a^{c^{-m}}\right]$ fixes $\Omega\setminus \bigcup_{n\in 12\Z}\Omega_n$
  pointwise for all $n\in \Z$ and $m \in \N$.
\end{lem}
\begin{proof}
Let $m \in \N$ be fixed. Using the definitions of $a, c$ and $d$ it is not difficult to check that for all $n \in \Z$ and $i,j \in \Z^{*}$, both $a^{c^m d}$ and $a^{c^{-m}}$ map the set 
$$A_{n,i,j}=\{n\}\times \{i\} \times \{j\} \times \Q^{*}$$
 to itself. Hence it suffices to consider the action of   $a^{c^m d}$ and $a^{c^{-m}}$ on any given $A_{n,i,j}$.
If $n \not \in 12\Z$, then $A_{n,i,j}$  is fixed pointwise by $a,c,$ and $d$ and we are done. So we may assume that $n\in 12\Z$. For simplicity, we will in fact assume that $n \in 24\Z$. The proof in the case that $n \in 24\Z + 12$ is analogous and omitted. 

If $(i,j) \not = (0,0)$, then, using \eqref{conjugate_support},
\begin{align*}
\supp\left(a^{c^{md}}\right)\cap A_{n,i,j}
&=\left(\supp(a) \cap A_{n,i,j}{(c^{md})}^{-1}\right)c^{m}d\\
&=\left(\supp(a) \cap A_{n,i,j - m}\right)c^md\\
&\subseteq \left(\{n\}\times \{i\}\times \{j - m\} \times (-1,1)\right)c^{m}d\\
&\subseteq \{n\}\times \{i\}\times \{j\} \times (1,3)
\end{align*}
and
\begin{align*}
\supp\left(a^{c^{-m}}\right)\cap A_{n,i,j}
&=\left(\supp(a) \cap A_{n,i,j}c^{m}\right)c^{-m}\\
&=\left(\supp(a) \cap A_{n,i,j + m}\right)c^{-m}\\
&\subseteq \left(\{n\}\times \{i\}\times \{j + m\} \times (-1,1)\right)c^{-m}\\
&\subseteq \{n\}\times \{i\}\times \{j\} \times (-1,1).
\end{align*}
Thus the supports of $a^{c^m d}$ and $a^{c^{-m}}$ are disjoint on $A_{n,i,j}$ if $(i,j)\not = (0,0)$, and so $\left[a^{c^md}, a^{c^{-m}}\right]$ fixes such $A_{n,i,j}$ pointwise.

It only remains to show that  $\left[a^{c^md}, a^{c^{-m}}\right]$ fixes $A_{n,0,0}\setminus \Omega_{n}$, i.e. all points of the form $(n,0,0,x)$ where $x \not \in (-1,1)$. In fact, using the definitions of $a,c$ and $d$ it is easy to verify that such $(n,0,0,x)$ are fixed under both  $a^{c^m d}$ and $a^{c^{-m}}$.
\end{proof}

We have shown that the sequence over the alphabet $\{f, f^{-1}, g, g^{-1}\}$ with $n$th term equal to
$$w_n=[(g\cdot g^{f^{-12}})^{(g^{f^{-4}})^ng^{f^{-28}}},
(g\cdot g^{f^{-12}})^{(g^{f^{-4}})^{-n}}]$$
is universal for $\stab(I_{48})$ as a subset of $\autq$.
Using \eqref{g_inverse} and \eqref{f_inverse} we can write
\begin{align*}
&g=f^{48}(f^{-48}g),\\ 
&g^{-1}=(f^{-48}g)f^{48}(f^{-48}g)f^{96}(f^{-48}g),\\
&f^{-1}=f^{47}(f^{-48}g)(f^{-48}g)f^{48}(f^{-48}g)f^{96}(f^{-48}g)
\end{align*}
and substituting these values in to $w_n$ yields a universal sequence over $\{f, f^{-48}g\}$ for
$\stab(I_{48})$ as a subset of $\autq$.  Combining this with equation
\eqref{eq-univ-sequence}, at the end of previous section, it is possible to
obtain an explicit universal sequence for $\autq$ over the alphabet $\{f,
f^{-48}g\}$. However, the resulting expression is too long to include here. 

%%%%%%%%%%%%%%%%%%%%%%%%%%%%

%\begin{acknowledgements}
%We would like to thank the anonymous referee, whose suggestions have significantly improved the presentation of this paper.
%\end{acknowledgements}

\bibliography{autq}{}
\bibliographystyle{plain}

\end{document}